\documentclass[12pt]{article}
\usepackage{a4, amssymb, amsmath, exscale}
\usepackage{makeidx, multicol}

\usepackage[mathscr]{eucal}
\usepackage[usenames,dvips]{color}
\usepackage{pstricks}

\newenvironment{evlist}[2]{
\begin{list}{}{
\setlength{\topsep}{0.5ex plus0.2ex minus0.1ex} 
\setlength{\leftmargin}{#1}
\setlength{\itemsep}{#2 plus0.2ex}
\setlength{\listparindent}{0pt}
\setlength{\parsep}{0ex plus0.2ex} }}
{\end{list}}

\newcommand{\Real}{\mathbb{R}}

\newcommand{\id}{\mathrm{id}}

\newcommand{\half}{{\textstyle\frac{1}{2}}}

\newcommand{\ol}[1]{\overline{#1}}
\newcommand{\co}[1]{{#1}^\diamond}

\newcommand{\interior}[1]{\mathrm{int}(#1)}

\renewcommand{\bmod}[1]{{#1}_{\setminus S}}

\newcommand{\inv}[2]{\Delta_{#2}(#1)}
\newcommand{\pinv}[2]{\Delta^{\!+}_{#2}(#1)}
\newcommand{\ninv}[2]{\Delta^{\!-}_{#2}(#1)}
\newcommand{\coinv}[2]{\Delta^{\hspace{-0.05em}\diamond}_{#2}(#1)}
\newcommand{\uninv}[2]{\underline{\Delta}^{\!-}_{#2}(#1)}
\newcommand{\gy}{\bullet}

\newcommand{\Minimal}[1]{{M}_{#1}}
\newcommand{\coMinimal}[1]{{M}^\diamond_{#1}}
\newcommand{\nMinimal}[1]{{N}_{#1}}
\newcommand{\conMinimal}[1]{{N}^\diamond_{#1}}
\newcommand{\Germ}[2]{G_{#2}(#1)}
\newcommand{\coGerm}[2]{G^{\hspace{0.1em}\diamond}_{#2}(#1)}

\newcommand{\Zed}[1]{\mathrm{Z}_{#1}}
\newcommand{\Inv}[1]{\mathcal{I}_{#1}}
\newcommand{\coInv}[1]{\mathcal{I}^\diamond_{#1}}
\newcommand{\fInv}[1]{\mathcal{D}_{#1}}
\newcommand{\cofInv}[1]{\mathcal{D}^{\hspace{0.05em}\diamond}_{#1}}

\newcommand{\definition}[1]{\textit{#1}}

\newcommand{\proof}{{\textit{Proof}\enspace}}

\newcommand{\eop}{\ \vbox{\hrule
                       \hbox{\vrule
                             \hskip 6pt
                             \vrule height 6pt width 0pt
                             \vrule}%
                       \hrule}%
                     \vspace{\medskipamount}
                }

\newtheorem{lemma}{Lemma}
\newtheorem{proposition}{Proposition}
\newtheorem{theorem}{Theorem}

\begin{document}

\title{Iterates of mappings which are almost\\ continuous and open}
\author{Chris Preston}
\date{\small{}}
\maketitle

\begin{quote}
This note presents an approach to studying the iterates of a mapping whose restriction to the complement of a 
finite set is continuous and open. The main examples to which the approach can be applied are
piecewise monotone mappings defined on an interval or a finite graph.
\end{quote}

\thispagestyle{empty}

If $X$ is a set and $g : X \to X$ a mapping then for each $n \ge 0$ the \definition{$n$-th iterate of $g$} will 
be denoted by $g^n$, i.e., $g^n : X \to X$, $n \ge 0$, are the mappings defined inductively by $g^0 = \id_X$ and
$g^n = g \circ g^{n-1}$ for each $n \ge 1$. Such a mapping $g$ will be thought of as describing a discrete dynamical
system, and in this interpretation the sequence $\{g^n(x)\}_{n \ge 0}$ is the 
\definition{orbit of the point $x \in X$ under $g$}.

We work here with a class of mappings which are defined as follows: Let $X$ be a topological space and 
$g : X \to X$ be a mapping. If $O$ is an open subset of $X$ then $g$ will be called \definition{regular on $O$}
if the restriction of $g$ to $O$, considered as a mapping from $O$ to $X$, is  both continuous and open ($O$ being 
endowed with the subspace topology). Hence, since $O$ is open, the requirement is that $g^{-1}(U) \cap O$ and
$g(U \cap O)$ must both be open subsets of $X$ for each open $U \subset X$. If $g$ is regular on each set in a 
family of open sets then it is clearly  also regular on their union. There is thus a largest open set $\Gamma_g$ on 
which $g$ is regular, and we say that $g$ is \definition{almost regular} if this set $\Gamma_g$ is dense, i.e., 
if its closure is the whole of $X$.

Piecewise monotone mappings of an interval are almost regular:
Let  $a ,\, b \in \Real$  with  $a < b$  and put  $I = [a,b]$. 
A mapping  $h : I \to I$  is said to be \definition{piecewise monotone} if there exists
$p \ge 0$  and  $a = d_0 < d_1 < \cdots < d_p < d_{p+1} = b$  such that  
$h$  is continuous and strictly monotone on each of the open intervals  $(d_k,d_{k+1})$, $k = 0,\,\ldots,\, p$.
Then $h$ is clearly regular on each of these intervals and hence it is regular on their union. But the complement 
of the union is the finite set $S_h = \{d_0,\ldots,d_{p+1}\}$ and thus
$h$ is almost regular with $\Gamma_h \supset I \setminus S_h$.
Note that $h$ is not assumed to be continuous at the points $d_0,\,d_1,\,\ldots,\,d_{p+1}$,
although the continuous case is much simpler to deal with.

There is a vast literature on such mappings. The topic in which we are interested (the asymptotic behaviour 
of `typical' orbits) is dealt with, for example, in Collet and Eckmann \cite{colleteckmann80},
Guckenheimer \cite{guckenheimer79},
Hofbauer \cite{hofbauer81} and \cite{hofbauer86}, Preston \cite{preston83} and \cite{preston88} and
Willms \cite{willms87}.

\begin{figure}
\begin{center}
\begin{pspicture}(0,0)(8,8)

\psset{xunit=0.8cm}
\psset{yunit=0.8cm}

\psframe(0.5,0.5)(9.5,9.5)

\psbezier(0.5,2.0)(1.5,6.0)(1.5,9.0)(2.0,9.0)
\psbezier(3.5,6.0)(2.5,6.0)(2.3,9.0)(2.0,9.0)
\psbezier(4.25,4.0)(4.0,3.0)(3.75,1.0)(3.5,1.0)
\psbezier(4.25,4.0)(4.5,6.0)(4.7,7.0)(5.0,7.0)
\psbezier(6.0,6.0)(5.3,6.5)(5.2,7.0)(5.0,7.0)
\psbezier(6.75,6.5)(6.3,6.3)(6.2, 5.0)(6.0,5.0)
\psbezier(6.75,6.5)(7.2,7.0)(7.3,8.0)(7.5,8.0)
\psbezier(9.5,8.0)(8.0,7.0)(7.8,7.0)(7.5,2.0)

\psline[linestyle=dashed](2.0,0.5)(2.0,9.0)
\psline[linestyle=dashed](3.5,0.5)(3.5,6.0)
\psline[linestyle=dashed](5.0,0.5)(5.0,7.0)
\psline[linestyle=dashed](6.0,0.5)(6.0,6.0)
\psline[linestyle=dashed](7.5,0.5)(7.5,8.0)

\uput{6pt}[d](0.5,0.4){$d_0$}
\uput{6pt}[d](2.0,0.4){$d_1$}
\uput{6pt}[d](3.5,0.4){$d_2$}
\uput{6pt}[d](7.5,0.4){$d_p$}
\uput{6pt}[d](9.5,0.4){$d_{p+1}$}

\end{pspicture}
\end{center}
\end{figure}

A generalisation of interval mappings are mappings defined on a graph with finitely many edges.
Here is an example: For $i = 1,\,2,\,3$ let $I_i = \{i\} \times [0,1]$ and let
$G$ be the quotient space $I_1 \cup I_2 \cup I_3/{\sim}$, where there are only two non-trivial 
equivalence classes given by $(1,0) \sim (2,0) \sim (3,0)$ and $(1,1) \sim (2,1) \sim (3,1)$. 

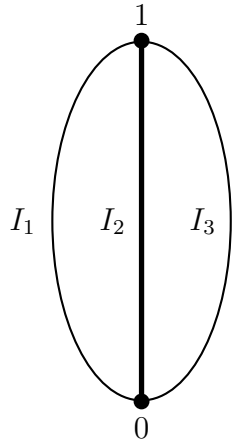
\begin{center}
\begin{pspicture}(0,0)(6,5)

\psset{xunit=0.6cm}
\psset{yunit=0.6cm}

\psline[linewidth=2pt](4.0,0.0)(4.0,8.0)
\psellipse(4.0,4.0)(2.0,4.0)

\psdot[dotsize=6pt](4.0,0.0)
\psdot[dotsize=6pt](4.0,8.0)

\uput{6pt}[l](2.0,4.0){$I_1$}
\uput{6pt}[l](4.0,4.0){$I_2$}
\uput{6pt}[l](6.0,4.0){$I_3$}
\uput{6pt}[d](4.0,0.0){$0$}
\uput{6pt}[u](4.0,8.0){$1$}

\end{pspicture}
\end{center}

Define a mapping $h : G \to G$ by
\[     h((a,x) = \left\{ \begin{array}{cl}
                           (a,2x) & \mbox{if $x \in [0,\half]$}\;,\\
                           (\sigma(a),2-2x) & \mbox{if $x \in [\half,1]$}
                           \end{array} \right.
\]
for $a = 1,\,2,\,3$, where $\sigma(1) = 2$, $\sigma(2) = 3$ and $\sigma(3) = 1$.
Then $h$ is almost regular with $\Gamma_h = G \setminus \{(1,\half),(2,\half),(3,\half)\}$. Moreover, $h$ is
continuous. This mapping $h$ is more-or-less equivalent to the piecewise monotone mapping
$h' : I \to I$, where $I = [0,3]$ and
\[     h'(x) = \left\{ \begin{array}{cl}
                           2x - a& \mbox{if $x \in (a,a + \half)$, $a = 0,\,1,\,2$}\;,\\
                           3(a + 1) - 2x& \mbox{if $x \in (a + \half, a + 1)$, $a = 0,\,1$}\;,\\
                             6  - 2x& \mbox{if $x \in (2\half, 3)$}\;,
                           \end{array} \right.
\]
and where $h'$ is assigned arbitrary values at the arguments $0,\,\half,\,1,\,1\half,\,2,\,2\half,\,3$.
For $k = 1,\,2,\,3$ the interval $[k-1,k] \subset I$ corresponds to the edge $I_k$ in $G$.

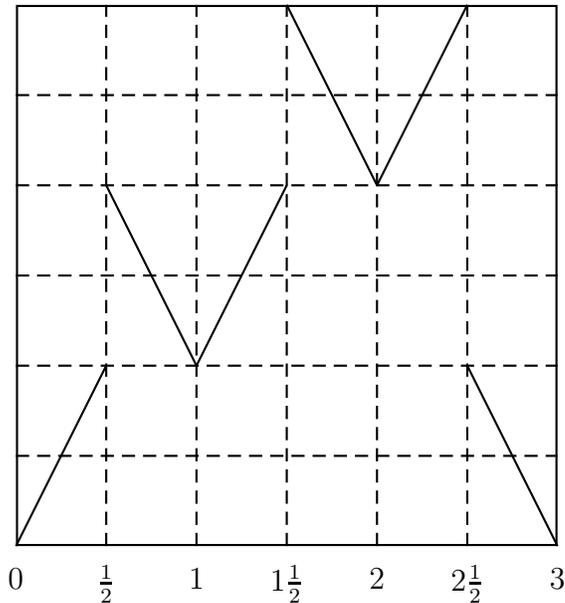
\begin{center}
\begin{pspicture}(0,0)(8,8)

\psset{xunit=0.8cm}
\psset{yunit=0.8cm}

\psframe(0.5,0.5)(9.5,9.5)

\psline(0.5,0.5)(2.0,3.5)
\psline(2.0,6.5)(3.5,3.5)
\psline(3.5,3.5)(5.0,6.5)
\psline(5.0,9.5)(6.5,6.5)
\psline(6.5,6.5)(8.0,9.5)
\psline(8.0,3.5)(9.5,0.5)

\psline[linestyle=dashed](0.5,2.0)(9.5,2.0)
\psline[linestyle=dashed](0.5,3.5)(9.5,3.5)
\psline[linestyle=dashed](0.5,5.0)(9.5,5.0)
\psline[linestyle=dashed](0.5,6.5)(9.5,6.5)
\psline[linestyle=dashed](0.5,8.0)(9.5,8.0)

\psline[linestyle=dashed](2.0,0.5)(2.0,9.5)
\psline[linestyle=dashed](3.5,0.5)(3.5,9.5)
\psline[linestyle=dashed](5.0,0.5)(5.0,9.5)
\psline[linestyle=dashed](6.5,0.5)(6.5,9.5)
\psline[linestyle=dashed](8.0,0.5)(8.0,9.5)

\uput{6pt}[d](0.5,0.4){$0$}
\uput{6pt}[d](2.0,0.45){$\half$}
\uput{6pt}[d](3.5,0.4){$1$}
\uput{6pt}[d](5.0,0.45){$1\half$}
\uput{6pt}[d](6.5,0.4){$2$}
\uput{6pt}[d](8.0,0.45){$2\half$}
\uput{6pt}[d](9.5,0.4){$3$}

\end{pspicture}
\end{center}

However, when working with the mapping $h'$ it is not easy to make use of the additional information of $h$ being 
continuous.

There is a now also a large literature on piecewise monotone mappings defined on finite graphs, see,
for example, Alseada, Llibre and Misiurewicz \cite{alm89}, Barge and Diamond \cite{bargediamond94} or
Llibre and Misiurewicz \cite{lm94}.

Each non-constant rational mapping $g : \Sigma \to \Sigma$ of the Riemann sphere $\Sigma$ into (and thus onto)
itself is regular, i.e., it is almost regular with $\Gamma_g = \Sigma$. However, the framework to be developed 
below does not really have much to say about this case.

Given an almost regular mapping $g : X \to X$, we are interested in describing the asymptotic behaviour of the 
orbit of a `typical' point  $x$ under $g$. In what follows $X$ will always be a complete metric space, and 
`typical' will be taken in the sense associated with the Baire category theorem.  Let us thus review the concepts 
involved here. (For more information Oxtoby \cite{oxboby71} is to be recommended.) The interior of a subset $A$ 
of $X$ will be denoted by $\interior{A}$, its closure by $\ol{A}$ and its boundary (i.e., the subset 
$\ol{A} \setminus \interior{A}$) by $\partial A$. Moreover, $A$ is \definition{dense} if $\ol{A} = X$ and  
\definition{nowhere dense} if $\interior{\ol{A}} = \varnothing$. In particular the boundary $\partial A$ of an 
open or a closed set $A$ is nowhere dense. The set $A$ is said to be \definition{meagre} if it can written as a 
countable union of nowhere dense sets. A set whose complement is meagre is called a \definition{residual set}; 
this is the case if and only if it contains a dense $G_\delta$-set, where a $G_\delta$-set is one which can be 
written as a countable intersection of open subsets of $X$. The Baire category theory (in its version for complete 
metric spaces) states that a countable intersection of dense open subsets of $X$ is itself dense, which implies
that every residual subset of $X$ is dense. Residual sets are clearly closed under taking countable intersections; 
they are also closed under taking finite unions.

A statement is considered to hold for `typical' points in $X$ if it holds on a residual set. For a given almost 
regular mapping $g : X \to X$ we thus want to make statements about the asymptotic behaviour of the orbit 
$\{g^n(x)\}_{n\ge 0}$ which hold for all points $x$ lying in some residual subset of $X$.

In these notes we present an approach for dealing with the iterates of an almost regular mapping $g$. Most of the 
results hold without any further assumptions. However, in one of the main results (Theorem~\ref{theorem_3}) 
$X \setminus \Gamma_g$ is required to be a finite set (which is the case for piecewise monotone mappings).

We also assume that $X$ is separable, and so there is a countable base for the topology, and that $X$ is perfect 
(meaning that there are no isolated points). This ensures that each finite subset of $X$ is nowhere dense.

\begin{lemma}~\label{lemma_1}
If $B \subset X$ is nowhere dense then so is $g^{-1}(B)$.
\end{lemma}

\proof 
Put $U = \Gamma_g$ and first consider a closed subset $E$ of $X$. Then $g^{-1}(E) \cap U$ is closed in the 
subspace topology on $U$ and so there exists a closed subset $F$ of $X$ such that $g^{-1}(E) \cap U = F \cap U$; 
in particular $g^{-1}(E) \subset F \cup (X \setminus U)$ and thus $\ol{g^{-1}(E)} \subset F \cup (X \setminus U)$.
Suppose now there exists a non-empty open set $O$ with $O \subset \ol{g^{-1}(E)}$. Then
$O \cap U \subset (F \cup (X \setminus U)) \cap U =  F \cap U = g^{-1}(E) \cap U \subset g^{-1}(E)$
and therefore $g(O \cap U) \subset E$. But $g(O \cap U)$ is open and $O \cap U \ne \varnothing$ (since 
$X \setminus U$ is nowhere dense), and hence $\interior{E} \ne \varnothing$. This shows that if 
$\interior{E} = \varnothing$ then $\interior{\,\ol{g^{-1}(E)}\,} = \varnothing$. It follows that if $B$ is nowhere 
dense, i.e., $\interior{\ol{B}} = \varnothing$ then $\interior{\,\ol{g^{-1}(\ol{B})}\,} = \varnothing$, which 
implies that $\interior{\,\ol{g^{-1}(B)}\,} = \varnothing$, i.e., that $g^{-1}(B)$ is nowhere dense.
\eop

For each $A \subset X$ put $\Lambda_g(A) = \{ x \in X : \mbox{$g^n(x) \notin A$ for all $n \ge 0$} \}$; then
$\Lambda_g(A)$ is $g$-invariant, meaning that $g(\Lambda_g(A)) \subset \Lambda_g(A)$, and thus the orbit 
$\{g^n(x)\}_{n\ge 0}$ of each point $x \in \Lambda_g(A)$ remains in $\Lambda_g(A)$. Moreover, if $S \subset X$ is 
nowhere dense then by Lemma~\ref{lemma_1} the set $X \setminus \Lambda_g(S) = \bigcup_{n\ge 0} g^{-n}(S)$ is 
meagre, and so $\Lambda_g(S)$ is a residual set.) Therefore, since we are only interested in the asymptotic 
behaviour of the orbits of  `typical' points in $X$, we can choose to ignore what happens with the orbits of the 
points in $X \setminus \Lambda_g(S)$. In particular, this means that it doesn't matter if we change the mapping 
$g$ on a nowhere dense set $S$, since this set is avoided by the orbits of all the points in the residual set 
$\Lambda_g(S)$.

We choose to modify the mapping $g$ as follows: Fix a closed nowhere dense set $S$ with 
$X \setminus \Gamma_g \subset S$, let $\gy$ be some element not in  $X$, put  $X_\gy = X \cup \{\gy\}$ and 
define a mapping $f : X_\gy \to X_\gy$ by 
\[
f(x) = \left\{ \begin{array}{cl}
                     g(x)  &\ \mbox{if $x \in X \setminus S$}\;,\\
                      \gy &\ \mbox{if $x \in S \cup \{\gy\}$}\;.
\end{array} \right.
\]
Thus $\Lambda_f(S) = \{ x \in X : \mbox{$f^n(x) \notin S$ for all $n \ge 0$} \} = \Lambda_g(S)$, by
Lemma~\ref{lemma_1} $\Lambda_f(S)$ is a residual set, and for each $x \in \Lambda_f(S)$ the orbit of $x$ 
under $f$ is the same as its orbit  under $g$. In the applications $S$ will always be taken to be 
$X \setminus \Gamma_g$, although taking $S$ somewhat larger may sometimes give additional information.
It turns out to be much more convenient to work with $f$ than with the mapping $g$, and so this is what we will do.

The set $X_\gy$ will not be endowed with a topology, and statements of a topological nature always refer to the 
topology on $X$, and thus must be statements about subsets of $X$. For example, the statement that a set $O$ is 
open means that $O$ is an open subset of $X$ (and in particular includes  the assertion that $O \subset X$).  
Since $X \setminus S$ is open the statement that $O$ is an open subset of $X \setminus S$ means that $O$ is a 
subset of $X \setminus S$ which is an open subset of $X$.

From the definition of $f$ it follows that $f(X \setminus S) \subset X$ and that the restriction of $f$ to 
$X \setminus S$ is equal to the restriction of $g$ to $X \setminus S$. Thus the restriction of $f$ to 
$X \setminus S$, considered as a mapping from $X \setminus S$ to $X$, is both continuous and open (where
$X \setminus S$ is again given the subspace topology). More explicitly, this means:

\begin{evlist}{32pt}{6pt}
\item[($A_1$)]
$f^{-1}(O)$ is an open subset of $X \setminus S$ for each open set $O$. 

\item[($A_2$)]
$f(O \setminus S)$ is open for each open set $O$. 
\end{evlist}

Condition ($A_1$) holds because $f^{-1}(A) =  g^{-1}(A) \setminus S$ for all $A \subset X$, and ($A_2$) holds 
because $f(A \setminus S) = g(A \setminus S)$ for all $A \subset X$. 

Since $f^{-1}(X) \subset X \setminus S$ it follows by induction that $f^{-n}(X) \subset X \setminus S$ for all 
$n \ge 1$. It then also follows from ($A_1$) by induction that if $O$ is open then $f^{-n}(O)$ is an open subset 
of $X \setminus S$ for each $n \ge 1$.

For each $A \subset X$ we denote the set $f(A \setminus S)$ by $\bmod{f}(A)$ (although some care is needed here, 
since there is no mapping $f_{\setminus S}$ involved). In particular ($A_2$) implies that $\bmod{f}(O)$ is an 
open subset of $X \setminus S$ for each open set $O$.

A subset $A$ of $X$ will be called \definition{$\bmod{f}$-invariant} if $\bmod{f}(A) \subset A$. Being more 
explicit, this means that $f(A \setminus S) \subset A$ must hold, and note that $f(A \setminus S) = X \cap f(A)$.
The subset $A$ is $\bmod{f}$-invariant if and only if $A \cup \{\gy\}$ is $f$-invariant, i.e., if and only if 
$f(A \cup \{\gy\}) \subset A \cup \{\gy\}$. This holds because if $f(A \cup \{\gy\}) \subset A \cup \{\gy\}$ then
\[ f(A\setminus S) = X \cap f(A) \subset X \cap f(A \cup \{\gy\}) \subset X \cap (A \cup \{\gy\}) = A\;,\]
and if $f(A \setminus S) \subset A$ then $f(A \cup \{\gy\}) = f(A \setminus S) \cup \{\gy\} \subset A \cup \{\gy\}$. 

The set $\Lambda_f(S) = \{ x \in X : \mbox{$f^n(x) \notin S$ for all $n \ge 0$} \}$ is $\bmod{f}$-invariant and 
for each $A \subset X$ the set $A \cap \Lambda_f(S)$ is $\bmod{f}$-invariant if and only if it is $g$-invariant.

The sets $\varnothing$ and $X$ are clearly $\bmod{f}$-invariant and arbitrary unions and intersections of 
$\bmod{f}$-invariant sets are again $\bmod{f}$-invariant.

\begin{lemma}\label{lemma_2}
If $A$ is $\bmod{f}$-invariant then $X \cap f^n(A) \subset A$ for all $n \ge 1$. Moreover, if $O$ is open
then so is $X \cap f^n(O)$. 
\end{lemma}

\proof
Since $X \cap f(B) = X \cap f(X \cap B)$  for all $B \subset X$ we have
\[
 X \cap f^{n+1}(A) = X \cap f(f^n(A)) = X \cap f(X \cap f^n(A)) = \bmod{f}(X \cap f^n(A)) 
\]
for all $A \subset X$ and all $n \ge 0$. If $A$ is invariant then $X \cap f(A) = \bmod{f}(A) \subset A$ and so  
it follows from the above and by induction on $n$ that $X \cap f^n(A) \subset A$ for all $n \ge 1$.
The final statement also follows by induction since $\bmod{f}(O')$ is open for each open set $O'$.
\eop

\begin{lemma}\label{lemma_3}
For each subset $A$ of $X$ the complement $X \setminus A$ is $\bmod{f}$-invariant if and only if 
$f^{-1}(A) \subset A$.
\end{lemma}

\proof 
The set $X \setminus A$ is $\bmod{f}$-invariant if and only if 
$f((X \setminus A) \setminus S) \subset X \setminus A$ and, since 
$(X \setminus A) \setminus S = X \setminus (A \cup S)$, this holds if and only if $f(x) \in X \setminus A$ whenever 
$x \in X \setminus (A \cup S)$, and thus if and only if $x \in A \cup S$ whenever $f(x) \in A$. Hence
$X \setminus A$ is $\bmod{f}$-invariant if and only if $f^{-1}(A) \subset A \cup S$. But 
$f^{-1}(A) \subset A \cup S$ holds if and only if $f^{-1}(A) \subset A$ does, since 
$f^{-1}(A) \subset X \setminus S$. 
\eop

If $X \setminus A$ is $\bmod{f}$-invariant then it follows from  Lemma~\ref{lemma_3} and by induction on $n$ 
that $f^{-n}(A) \subset A$ for all $n \ge 1$, since $f^{-(n+1)}(A) = f^{-1}(f^{-n}(A))$ for all $n \ge 0$.

The following lemma gives the two properties of $\bmod{f}$-invariant sets which play a fundamental role in the 
analysis of the iterates of $f$.

\begin{lemma}\label{lemma_4}  
If a subset $A$ of $X$ is $\bmod{f}$-invariant then so are its interior $\interior{A}$ and its closure $\ol{A}$.
\end{lemma}

\proof
Let $A$ be $\bmod{f}$-invariant; then $\bmod{f}(\interior{A}) \subset \bmod{f}(A) \subset  A$ and by ($A_2$) 
the set $\bmod{f}(\interior{A})$ is open. Hence $\bmod{f}(\interior{A}) \subset  \interior{A}$, i.e., 
$\interior{A}$ is $\bmod{f}$-invariant. Now consider $x \in \ol{A} \setminus S$, and let $O$ be any open set 
containing $f(x)$. Then $x$ lies in the set $f^{-1}(O)$, which by ($A_1$) is open. Hence 
$f^{-1}(O) \cap A \ne \varnothing$ (since $x \in \ol{A}$), and so there exists $y \in f^{-1}(O) \cap A$. In 
particular $y \in A \setminus S$ (since $f(y) \in O \subset X$), and it follows that $f(y) \in A$, since $A$ is 
$\bmod{f}$-invariant, which implies that $f(y) \in O \cap A$. Therefore $A$ intersects each open set containing 
$f(x)$, and thus $f(x) \in \ol{A}$. This shows that $\bmod{f}(\ol{A}) = f(\ol{A} \setminus S) \subset \ol{A}$, 
and hence that $\ol{A}$ is $\bmod{f}$-invariant.
\eop

A subset $A \subset X$ will be called \definition{fully $\bmod{f}$-invariant} if both $A$ and its complement 
$X \setminus A$ are $\bmod{f}$-invariant. The sets $\varnothing$ and $X$ are fully $\bmod{f}$-invariant and
arbitrary unions and intersections of fully $\bmod{f}$-invariant sets are again  fully $\bmod{f}$-invariant. 
Moreover, by definition they are closed under taking complements. If $A$ is fully $\bmod{f}$-invariant then 
Lemma~\ref{lemma_4} shows that the interior $\interior{A}$ and the closure $\ol{A}$ of $A$ are also fully 
$\bmod{f}$-invariant, since $X \setminus \interior{A} = \ol{X \setminus A}$ and 
$X \setminus \ol{A} = \interior{X \setminus A}$.

The set of all fully $\bmod{f}$-invariant open subsets of $X$ will be denoted by $\fInv{f}$. Our aim is to find a 
decomposition of $X$ into disjoint elements $\{D_\alpha\}$ of $\fInv{f}$ such that the union 
$\bigcup_{\alpha} D_\alpha$ is dense in $X$ and such that the behaviour of the iterates of $f$ on each of the 
components $D_\alpha$ allows some kind of reasonable description. 

In fact, it is convenient to only consider the regular open sets in $\fInv{f}$: For each subset $B \subset X$ 
put $\co{B} = \interior{\ol{B}}$; a set $O$ is said to be \definition{regular open} if $O = \co{O}$.
(Such sets were introduced by Kuratowksi in \cite{kuratowski22} and have proved useful in various settings.) 
In particular $B$ is nowhere dense if and only if $\co{B} = \varnothing$. If $A \subset B$ then clearly 
$\co{A} \subset \co{B}$. The elementary properties of the $\diamond$-operation which we will make use of are 
listed in the following lemma:

\begin{lemma}\label{lemma_5}
(1)\enskip 
$\co{(\co{B})} = \co{B}$ for all $B \subset X$. 

(2)\enskip 
A set $O$ is regular open if and only if
$O = \co{B}$ for some $B$, which is the case  if and only if $O = \interior{F}$ for some closed set $F$.

(3)\enskip
If $O_1,\,O_2$ are regular open then so is $O_1 \cap O_2$.

(4)\enskip
If $O$ is open (and in particular if $O$ is regular open) then $X \setminus \ol{O}$ is regular open.

(5)\enskip
If $O_1,\,O_2$ are regular open then so is $O_2 \setminus \ol{O_1}$. Moreover, if $O_1$ is a proper subset of
$O_2$ then $O_2 \setminus \ol{O_1}$ is non-empty.

(6)\enskip
If $O \subset X$ is open then $O \subset \co{O}$ and $\co{O} \setminus O \subset \partial O$ (the boundary of $O$), 
and in particular $\co{O} \setminus O$ is nowhere dense.

(7)\enskip
If $O_1,\,O_2$ are open then $O_1 \cap O_2 \ne \varnothing$ if and only if $O_1 \cap \co{O}_2 \ne \varnothing$ 
which is the case if and only if $\co{O}_1 \cap \co{O}_2 \ne \varnothing$.
\end{lemma}

\proof 
(1)\enskip 
$\co{B} = \interior{\ol{B}} 
= \interior{\interior{\ol{B}}} \subset \interior{\ol{\interior{\ol{B}}}} = \co{(\co{B})}$,
and on the other hand
$\co{(\co{B})} = \interior{\ol{\interior{\ol{B}}}} \subset \interior{\ol{\ol{B}}} = \interior{\ol{B}} = \co{B}$. 

(2)\enskip
If $O = \co{B}$ for some $B$ then $\co{O} = \co{(\co{B})} = \co{B} = O$ and so $O$ is regular open. Conversely, 
if $O$ is regular open then $O = \co{B}$ with $B = O$. The rest follows since 
$\interior{F} = \interior{\ol{F}} = \co{F}$ whenever $F$ is open.

(3)\enskip
We have $\co{(O_1 \cap O_2)} \subset \co{O}_1 = O_1$, since $O_1 \cap O_2 \subset O_1$, and in the same way
$\co{(O_1 \cap O_2)} \subset  O_2$. Thus $\co{(O_1 \cap O_2)} \subset  O_1 \cap O_2$. But $O \subset \co{O}$
for each open set $O$ and so also $O_1 \cap O_2 \subset \co{(O_1 \cap O_2)}$. Hence 
$O_1 \cap O_2 = \co{(O_1 \cap O_2)}$, i.e., $O_1 \cap O_2$ is regular open.

(4)\enskip
This follows from (2), since $X \setminus \ol{O} = \interior{X \setminus O}$ and $X \setminus O$ is closed.

(5)\enskip
By (3) and (4) $O_2 \setminus \ol{O_1} = O_2 \cap (X \setminus \ol{O_1})$ is regular open. Now if 
$O_2 \subset \ol{O_1}$ then $O_2 \subset \interior{\ol{O_1}} = \co{O}_1 = O_1$ and so $O_1$ is not a proper subset 
of $O_2$. Thus if $O_1$ is a proper subset of $O_2$ then $O_2 \setminus \ol{O_1}$ is non-empty.

(6)\enskip
Since $O \subset \ol{O}$ and $O$ is open it follows that $O \subset \interior{\ol{O}} = \co{O}$. Moreover,
$\co{O} \setminus O \subset \ol{O} \setminus O = \partial O$.

(7)\enskip
If $O_1 \cap O_2 = \varnothing$ then $O_1 \cap \ol{O_2} = \varnothing$ and hence 
$O_1 \cap \co{O}_2 = O_1 \cap \interior{\ol{O_2}} = \varnothing$. In the same way it then follows that 
$\co{O}_1 \cap \co{O}_2 = \varnothing$. The converse is clear, since $O_1 \subset \co{O}_2$ and 
$O_2 \subset \co{O}_2$.
\eop

We are working here with open sets where more typically closed sets would be used, and the
$\diamond$-operation can be thought of taking on the role normally played by the closure operation. 

Denote by $\cofInv{f}$ the set of non-empty regular open elements in $\fInv{f}$. Two important properties of 
elements $D_1,\,D_2 \in \cofInv{f}$ are the following: If $D_1 \cap D_2 \ne \varnothing$ then by 
Lemma~\ref{lemma_5}~(3) $D_1 \cap D_2 \in \cofInv{f}$. Moreover, if $D_1$ is a proper subset of $D_2$ then
by Lemma~\ref{lemma_5}~(5) $D_2 \setminus \ol{D_1}$ is non-empty, and thus an element of $\cofInv{f}$.

A set $D \in \cofInv{f}$ will be called \definition{minimal} if the only element of $\cofInv{f}$ which is a subset 
of $D$ is $D$ itself. Equivalently, $D \in \cofInv{f}$ is minimal if $D \subset E$ whenever $E \in \cofInv{f}$ 
with $E \cap D \ne \varnothing$ (since then $E \cap D \in \cofInv{f}$). This implies that if 
$D_1,\,D_2 \in \cofInv{f}$ are minimal with $D_1 \ne D_2$ then $D_1 \cap D_2 = \varnothing$, i.e., minimal elements 
are either equal or disjoint. In particular, the set of minimal elements in $\cofInv{f}$ is countable, since we 
are assuming $X$ is separable.

Denote by $\Minimal{f}$ the union of all the minimal elements in $\cofInv{f}$; thus $\Minimal{f} \in \fInv{f}$. Put 
$\nMinimal{f} = X \setminus \ol{\Minimal{f}}$ and so by Lemma~\ref{lemma_5}~(4) $\nMinimal{f}$ is regular 
open, which means that $\nMinimal{f}$ is either empty or an element of $\cofInv{f}$. Moreover, $\Minimal{f}$ and 
$\nMinimal{f}$ are disjoint and their union $\Minimal{f} \cup \nMinimal{f}$ is a dense open subset of $X$, since 
$X \setminus (\Minimal{f} \cup \nMinimal{f}) = \ol{\Minimal{f}} \setminus \Minimal{f}$ is the boundary of the open 
set $\Minimal{f}$.

For each $x \in X$ let $\cofInv{f}(x)$ be the set of all elements in  $\cofInv{f}$ which contain $x$, and put 
$\Germ{x}{f} = \bigcap_{D \in \cofInv{f}(x)} \ol{D}$, i.e., $\Germ{x}{f}$ is the intersection of the closures of the 
elements in  $\cofInv{f}(x)$. Hence $\Germ{x}{f}$ is a closed fully $\bmod{f}$-invariant set containing $x$ (since 
$X \in \cofInv{f}(x)$). We write $\coGerm{x}{f}$ instead of $\co{(\Germ{x}{f})}$; thus 
$\coGerm{x}{f} = \interior{\Germ{x}{f}}$, since $\Germ{x}{f}$ is closed, and $\coGerm{x}{f}$ is either empty or 
an element of $\cofInv{f}$.

Here is a more explicit description of the set $\Germ{x}{f}$: Since $X$ is fully $\bmod{f}$-invariant and
arbitrary intersections of fully $\bmod{f}$-invariant sets are fully $\bmod{f}$-invariant there is a least 
fully $\bmod{f}$-invariant set containing  a given subset $A$ of $X$, and this set will be denoted by 
$\inv{A}{f}$. If $O$ is open then  by Lemma~\ref{lemma_4} $\inv{O}{f}$ is also open, since then 
$\interior{\inv{O}{f}}$ is a fully $\bmod{f}$-invariant set containing $O$. If $D \in \cofInv{f}(x)$ then 
$B_{1/n}(x) \subset D$ for some $n \ge 1$ (with $B_r(x)$ the open ball of radius $r$ and centre $x$) and hence 
$\inv{B_{1/n}(x)}{f} \subset D$; on the other hand $\inv{B_{1/n}(x)}{f} \in \cofInv{f}(x)$ for each $n$. It follows 
that $\Germ{x}{f} = \bigcap_{n \ge 1} \ol{\inv{B_{1/n}(x)}{f}}$, which represents $\Germ{x}{f}$ as the intersection 
of a decreasing sequence of sets. For each $A \subset X$ it is convenient to write $\coinv{A}{f}$ instead of 
$\co{(\inv{A}{f})}$.

\begin{theorem}\label{theorem_1}
(1)\enskip
If $D$ is a minimal element of $\cofInv{f}$ then $\coGerm{x}{f} = D$ for each $x \in D$, and in particular
$x \in \coGerm{x}{f}$.

(2)\enskip
There exists a meagre subset $Z$ of $X$ such that $\coGerm{x}{f} = \varnothing$ (i.e., such that $\Germ{x}{f}$ is 
nowhere dense) for all $x \in \nMinimal{f} \setminus Z$.

(3)\enskip
$\Minimal{f} \subset \interior{\{ x \in X : \coGerm{x}{f} \ne \varnothing\}} \subset \coMinimal{f}$.
\end{theorem}

\proof
(1)\enskip
If $E \in \cofInv{f}(x)$ then $D \subset E$, since $D$ is minimal; also, $D \in \cofInv{f}(x)$. Thus 
$\Germ{x}{f} = \ol{D}$ and hence $\coGerm{x}{f} = \co{D} = D$.

(2)\enskip
Let $\{O_n\}_{n\ge 1}$ be a sequence of non-empty open sets such that each non-empty open set contains some $O_n$ 
(which exists since the topology on $X$ has a countable base). For each $n \ge 1$ put $D_n = \coinv{O_n}{f}$ and 
let $Z = \bigcup_{n \ge 1} (\ol{D_n} \setminus D_n)$; then $Z$ is a meagre subset of $X$, since 
$\ol{D_n} \setminus D_n$, as the boundary of the open set $D_n$, is nowhere dense for each $n$. Now let 
$x \in \nMinimal{f} \setminus Z$ and suppose $\coGerm{x}{f} \ne \varnothing$; then 
$\coGerm{x}{f} \subset \conMinimal{f} = \nMinimal{f}$, since $D \cap \nMinimal{f} \in \cofInv{f}(x)$ for each
$D \in \cofInv{f}(x)$. Hence $\coGerm{x}{f}$ is not minimal and so there exists $D \in \cofInv{f}$ with $D$ a proper 
subset of $\coGerm{x}{f}$. Let $n \ge 1$ be such that $O_n \subset D$; then $\inv{O_n}{f} \subset D$, since  
$\inv{O_n}{f}$ is the least fully $\bmod{f}$-invariant set containing $O_n$ and thus also 
$D_n = \coinv{O_n}{f} \subset \co{D} = D$. By Lemma~\ref{lemma_5}~(5) 
$\coGerm{x}{f} \setminus \ol{D} \ne \varnothing$, and 
$\coGerm{x}{f} \setminus \ol{D_n} \supset \coGerm{x}{f} \setminus \ol{D}$, which shows that
$\coGerm{x}{f} \setminus \ol{D_n} \ne \varnothing$. In particular $E_n = X \setminus \ol{D_n}$ is non-empty and so 
an element of $\cofInv{f}$. However, $x$ is either in $D_n$ or $E_n$, since 
$X \setminus (D_n \cup E_n) = \ol{D_n} \setminus D_n \subset Z$, which means that one of $D_n$ or $E_n$ is an 
element of $\cofInv{f}(x)$. Suppose $D_n \in \cofInv{f}(x)$; then 
$\coGerm{x}{f} \subset \interior{\ol{D_n}} = \co{D}_n = D_n$, and this is not possible, since $D_n$ is a proper 
subset of $\coGerm{x}{f}$. But if $E_n \in \cofInv{f}(x)$ then 
$\coGerm{x}{f} \subset \interior{\ol{E_n}} = \co{E}_n = E_n$, which again is not possible because 
$D_n \subset \coGerm{x}{f} \setminus E_n$. It therefore follows that $\coGerm{x}{f} = \varnothing$. 

(3)\enskip
Put $U = \interior{\{ x \in X : \coGerm{x}{f} \ne \varnothing\}}$. By (1) $\coGerm{x}{f} \ne \varnothing$ for all 
$x \in \Minimal{f}$ and so $\Minimal{f} \subset U$. Now let $x \in \nMinimal{f}$ and let $O$ be any open set 
containing $x$; then $O \cap \nMinimal{f} \ne \varnothing$. Thus 
$O \cap (\nMinimal{f} \setminus Z) = (O \cap \nMinimal{f}) \cap (X \setminus Z) \ne \varnothing$, since 
the residual set $X \setminus Z$ is dense. It follows that $\coGerm{y}{f} = \varnothing$ for some $y \in O$, which 
implies that $x \notin U$, and shows that $\nMinimal{f} \cap U = \varnothing$. Therefore
$U \subset \interior{X \setminus \nMinimal{f}} = X \setminus \ol{\nMinimal{f}}$. But 
$\nMinimal{f} = X \setminus \ol{\Minimal{f}}$ and so
$X \setminus \ol{\nMinimal{f}} = X \setminus \ol{X \setminus \ol{\Minimal{f}}}
= X \setminus (X \setminus \coMinimal{f}) = \coMinimal{f}$. \eop

If $x \in \partial D$ for some minimal element $D$ of $\cofInv{f}$ then $\coGerm{x}{f} \supset D$, which implies
in particular that $\coGerm{x}{f} \ne \varnothing$.

For each $x \in X$ the singleton set $\{x\}$ will be denoted just by $x$, so for example $\inv{x}{f}$ is the 
least fully $\bmod{f}$-invariant set containing $x$. It is easily checked that
\[
\inv{x}{f} = \{ y \in X : \mbox{there exist $m,\,n \ge 0$ with $f^m(y) \in X$ and $f^m(y) = f^n(x)$} \}\;.
\]
Note that if $A$ is any fully $\bmod{f}$-invariant set with $A \cap \inv{x}{f} \ne \varnothing$ then $x \in A$. 
(There exists $y \in A$ and $m,\,n \ge 0$ with $f^m(y) \in X$ and $f^m(y) = f^n(x)$ and thus by 
Lemma~\ref{lemma_2} $f^n(x) = f^m(y) \in A$; Lemma~\ref{lemma_3} then shows that $x \in A$.)

Theorem~\ref{theorem_1}~(2) implies in particular that there exists a meagre subset $Z$ of $X$ such that 
$\inv{x}{f}$ is nowhere dense for all $x \in \nMinimal{f} \setminus Z$ (since $\ol{\inv{x}{f}} \subset \Germ{x}{f}$ 
for all $x$). In fact we have the following:

\begin{proposition}\label{prop_1}
The set $\inv{x}{f}$ is nowhere dense for all $x \in \nMinimal{f}$.
\end{proposition}

\proof 
Let $x \in \nMinimal{f}$; then, since $\inv{x}{f} \subset \nMinimal{f}$ it follows that 
$\coinv{x}{f} \subset \conMinimal{f} = \nMinimal{f}$. Suppose $\coinv{x}{f} \ne \varnothing$, which means 
$\coinv{x}{f}$ is an element of $\cofInv{f}$ and a subset of $\nMinimal{f}$. Hence $\coinv{x}{f}$ is not minimal 
and so there exists $D \in \cofInv{f}$ with $D$ a proper subset of $\coinv{x}{f}$, and then by 
Lemma~\ref{lemma_5}~(5) $D' = \coinv{x}{f} \setminus \ol{D}$ is also an element of $\cofInv{f}$. It follows 
that each of the open sets $D$ and $D'$ intersects $\ol{\inv{x}{f}}$ and hence they each intersect $\inv{x}{f}$. 
But $D$ and $D'$ are both fully $\bmod{f}$-invariant, which implies they both contain $x$. This is not possible 
because they are disjoint, which shows that $\coinv{x}{f} = \varnothing$.
\eop

We next look at what goes on inside a minimal element of $\cofInv{f}$, and begin by considering  which sets should 
correspond to the minimal elements of $\cofInv{f}$ in the case of $\bmod{f}$-invariant (rather than fully 
$\bmod{f}$-invariant) sets.

Denote by $\Inv{f}$ the set of all $\bmod{f}$-invariant open subsets of $X$ and by $\coInv{f}$ the set of non-empty
regular open sets in $\Inv{f}$. Note that if $U_1,\,U_2 \in \coInv{f}$ with $U_1 \cap U_2 \ne \varnothing$ then by 
Lemma~\ref{lemma_5}~(3) $U_1 \cap U_2$ is again an element of $\coInv{f}$. An element $U \in \coInv{f}$ is 
said to be \definition{transitive} if the only $V \in \coInv{f}$ with $V \subset U$ is $U$ itself. Thus $U$ is 
transitive if and only if $U \subset V$ for each $V \in \coInv{f}$ with $V \cap U \ne \varnothing$ (since 
$V \cap U$ is then an element of $\coInv{f}$ with $V \cap U \subset U$). This implies in particular that 
transitive elements of $\coInv{f}$ are either equal or disjoint.

The transitive elements of $\coInv{f}$ are, of course, just the minimal elements of $\coInv{f}$. The denotation
`transitive' is employed to avoid the confusion of having two different usages for the term `minimal' and because
the property corresponds to the usual one of being topologically transitive (but expressed here for open rather than 
closed sets); see, for example, Walters \cite{walters82}.

A decreasing sequence $\{U_n\}_{n \ge 1}$ of sets from $\coInv{f}$ will be called an 
\definition{asymptotically transitive cascade}, or for short just a \definition{transitive cascade}, if for each 
$V \in \coInv{f}$ with $V \cap U_1 \ne \varnothing$ there exists $m \ge 1$ such that $U_m \subset V$. In 
particular, if $U \in \coInv{f}$ is transitive then the constant sequence $\{U_n\}_{n \ge 1}$ with $U = U_n$ for 
all $n \ge 1$ is a transitive cascade, and it will be denoted by $\{U\}$.
We will see that each minimal element of $\cofInv{f}$ contains an essentially unique transitive cascade and that, 
conversely, each transitive cascade is contained in a unique minimal element of $\cofInv{f}$.

For each transitive cascade $\gamma = \{U_n\}_{n \ge 1}$ we put $C_\gamma = \bigcap_{n \ge 1} \ol{U_n}$ and call 
$C_\gamma$ the \definition{core of $\gamma$}. If $X$ is compact then $C_\gamma$ is non-empty, but in general this 
need not be the case. The core $C_\gamma$ is closed and $\bmod{f}$-invariant, and hence 
$\co{C}_\gamma = \interior{C_\gamma}$ is either empty or an element of $\coInv{f}$. Moreover, 
$\co{C}_\gamma \subset U_n$ for all $n \ge 1$, since 
$\co{C}_\gamma = \interior{C_\gamma} \subset \interior{\ol{U_n}} = \co{U}_n = U_n$. If $U$ is a transitive element 
of $\coInv{f}$ then of course $C_{\{U\}} = \ol{U}$ and so $\co{C}_{\{U\}} = U$.

\begin{proposition}\label{prop_2}
If $\gamma = \{U_n\}_{n \ge }$ is a transitive cascade with $\co{C}_\gamma \ne \varnothing$ then $U_m$ is 
transitive for some $m \ge 1$ and $\co{C}_\gamma = U_m$. In particular, for each transitive cascade $\gamma$ the 
set $\co{C}_\gamma$ is either empty or a transitive element of $\coInv{f}$.
\end{proposition}

\proof 
We have $\co{C}_\gamma \in \coInv{f}$ and $\co{C}_\gamma \cap U_1 = \co{C}_\gamma \ne \varnothing$. Thus 
$U_m \subset \co{C}_\gamma$ for some $m \ge 1$, and thus $U_m \subset  \co{C}_\gamma \subset U_m$, i.e., 
$\co{C}_\gamma = U_m$. Now if $U$ is any element of $\coInv{f}$ with $U \cap U_m \ne \varnothing$ then 
$U_p \subset U$ for some $p \ge 1$ and hence $U_m = \co{C}_\gamma \subset U_p \subset U$. Therefore $U_m$ is 
transitive.
\eop

Let $\gamma = \{U_n\}_{n \ge }$ be a transitive cascade and consider the family $\mathcal{S}_\gamma$ consisting 
of those elements $U \in \coInv{f}$ which contain $U_1$ and are such that if $V \in \coInv{f}$ with 
$V \cap U \ne \varnothing$ then $U_m \subset V$ for some $m \ge 1$. (Thus $U$ being in $\mathcal{S}_\gamma$ means 
that $\{U'_n\}_{n \ge }$ is still a transitive cascade, where $U'_1 = U$ and $U'_{n+1} = U_n$ for all $n \ge 1$.) 
Now $\mathcal{S}_\gamma$ contains $U_1$ and thus is non-empty, and so by Lemma~\ref{lemma_5}~(7) the set 
$D_\gamma = \co{W}$, where $W$ is the union of all the sets in $\mathcal{S}_\gamma$, is again an element of 
$\mathcal{S}_\gamma$, and hence it is the largest element. We call $D_\gamma$ the \definition{domain of $\gamma$}.

If $U$ is a transitive element of $\coInv{f}$ then we write $D_U$ instead of $D_{\{U\}}$ and call $D_U$ the
\definition{domain of $U$}. Thus $D_U$ is the union of all the sets in the family $\mathcal{S}_U$ consisting 
of those elements $V \in \coInv{f}$ which contain $U$ and are such that if $V' \in \coInv{f}$ with 
$V' \cap V \ne \varnothing$ then $U \subset V'$.

\begin{theorem}\label{theorem_2}
The domain of a transitive cascade is a minimal element of $\cofInv{f}$. Conversely, each minimal element of 
$\cofInv{f}$ is the domain of some transitive cascade.
\end{theorem}

\proof 
This requires some preparation. For each subset $A \subset X$ put
\[  
\ninv{A}{f} = \{ x \in X : \mbox{$f^n(x) \in A$ for some $n \ge 0$} \}\;, 
\]
and hence $\ninv{A}{f} = \bigcup_{n \ge 0} f^{-n}(A)$, since $f^{-n}(A) \subset X$ for each $n$. In particular, if 
$O$ is open then so is $\ninv{O}{f}$.

\begin{lemma}\label{lemma_6}
(1)\enskip
The set $X \setminus \ninv{A}{f}$ is $\bmod{f}$-invariant for each $A \subset X$.

(2)\enskip
The set $X \setminus A$ is $\bmod{f}$-invariant if and only if $A = \ninv{A}{f}$. 

(3)\enskip
If $A$ is $\bmod{f}$-invariant then so is $\ninv{A}{f}$ and $\inv{A}{f} = \ninv{A}{f}$.

(4)\enskip
If $A$ and $B$ are $\bmod{f}$-invariant sets then $\inv{A}{f} \cap \inv{B}{f} \ne \varnothing$ if and only if 
$A \cap B \ne \varnothing$, which in turn holds if and only if  $A \cap \inv{B}{f} \ne \varnothing$.
\end{lemma}

\proof 
(1)\enskip
$f^{-1}(\ninv{A}{f}) = \{ x \in X : \mbox{$f^n(x) \in A$ for some $n \ge 1$} \} \subset \ninv{A}{f}$ and so by 
Lemma~\ref{lemma_2} $X \setminus \ninv{A}{f}$ is $\bmod{f}$-invariant. 

(2)\enskip
If $A = \ninv{A}{f}$ then by (1) $X \setminus A$ is $\bmod{f}$-invariant. Conversely, if  $X \setminus A$ is 
$\bmod{f}$-invariant then by Lemma~\ref{lemma_2} and induction $f^{-n}(A) \subset A$ for all $n \ge 1$ and 
hence $A \subset \ninv{A}{f} = \bigcup_{n \ge 0} f^{-n}(A)  \subset A$, i.e., $A = \ninv{A}{f}$. 

(3)\enskip
Let $x \in \ninv{A}{f} \setminus S$; then $f^n(x) \in A$ for some $n \ge 0$. If $x \in A$ then 
$x \in A \setminus S$ and so $f(x) \in A$, since $A$ is $\bmod{f}$-invariant. On the other hand, if $f^n(x) \in A$ 
with $n \ge 1$ then $f^m(f(x)) \in A$ with $m = n - 1$. In both cases $f(x) \in \ninv{A}{f}$ and so 
$\bmod{f}(\ninv{A}{f}) \subset \ninv{A}{f}$, i.e., $\ninv{A}{f}$ is $\bmod{f}$-invariant. Together with (1) this 
shows $\ninv{A}{f}$ is fully $\bmod{f}$-invariant and hence $\inv{A}{f} \subset \ninv{A}{f}$, since clearly 
$A \subset \ninv{A}{f}$. But by (2) $\ninv{\inv{A}{f}}{f} = \inv{A}{f}$, since $\inv{A}{f}$ is fully 
$\bmod{f}$-invariant, and $A \subset \inv{A}{f}$ and it follows that
$\ninv{A}{f} \subset \ninv{\inv{A}{f}}{f} = \inv{A}{f}$. Therefore $\inv{A}{f} = \ninv{A}{f}$.

(4)\enskip
Assume  $A \cap \inv{B}{f} \ne \varnothing$ and let $x \in A \cap \inv{B}{f}$. By (3) $x \in \ninv{B}{f}$ and so 
$f^n(x) \in B$ for some $n \ge 0$, thus $f^n(x) \in X \cap f^n(A)$ and hence by Lemma~\ref{lemma_2} 
$f^n(x) \in A$, i.e., $f^n(x) \in A \cap B$. Therefore $A \cap B \ne \varnothing$ whenever 
$A \cap \inv{B}{f} \ne \varnothing$. But $\inv{B}{f}$ is $\bmod{f}$-invariant, and applying this to the sets $A$ 
and $\inv{B}{f}$ then shows that $A \cap \inv{B}{f} \ne \varnothing$ whenever 
$\inv{A}{f} \cap \inv{B}{f} \ne \varnothing$. Finally, if $A \cap B \ne \varnothing$ then clearly both 
$A \cap \inv{B}{f} \ne \varnothing$ and $\inv{A}{f} \cap \inv{B}{f} \ne \varnothing$.
\eop

\begin{lemma}\label{lemma_7}
For each transitive cascade $\gamma$ the domain $D_\gamma$ is an element of $\cofInv{f}$. 
\end{lemma}

\proof
If $V \in \coInv{f}$ with $V \cap \coinv{D_\gamma}{f} \ne \varnothing$ then by Lemma~\ref{lemma_5}~(7) 
$V \cap \inv{D_\gamma}{f} \ne \varnothing$ and hence by Lemma~\ref{lemma_6}~(4) 
$V \cap D_\gamma \ne \varnothing$. It follows that $\coinv{D_\gamma}{f} \in \mathcal{S}_\gamma$, and thus that
$D_\gamma = \coinv{D_\gamma}{f}$ (since $D_\gamma \subset \coinv{D_\gamma}{f}$). But 
$\coinv{D_\gamma}{f} \in \cofInv{f}$ and so $D_\gamma \in \cofInv{f}$. 
\eop

\begin{lemma}\label{lemma_8}
If $\gamma = \{U_n\}_{n \ge 1}$ is a transitive cascade then $\coinv{U_n}{f} = D_\gamma$ for each $n \ge 1$.
\end{lemma}

\proof 
Consider $n$ to be fixed. By Lemma~\ref{lemma_7} we have $\coinv{U_n}{f} \subset \coinv{D_\gamma}{f} = D_\gamma$. 
Now put $U' = X \setminus \ol{\inv{U_n}{f}}$ (and so $U'$ is either empty or an element of $\coInv{f}$) and 
suppose the regular open set $W = D_\gamma \setminus \ol{\inv{U_n}{f}} = D_\gamma \cap U'$ is non-empty. Then we have 
$D_\gamma \cap U' \in \coInv{f}$, thus $U_m \subset D_\gamma \cap U'$ for some $m \ge 1$ and in particular 
$U_m \subset U'$. However, $U_n \subset \ol{\inv{U_n}{f}} = X \setminus U'$, which is not possible, since either 
$U_n \subset U_m$ or $U_m \subset U_n$. Hence $W = \varnothing$, i.e., $D_\gamma \subset \ol{\inv{U_n}{f}}$, which 
implies that $D_\gamma = \interior{D_\gamma} \subset \coinv{U_n}{f}$.
\eop

\begin{lemma}\label{lemma_9}
If $D \in \cofInv{f}$ is minimal then $U_1 \cap U_2 \ne \varnothing$ for all elements $U_1,\,U_2$ of $\Inv{f}$ 
contained in $D$. 
\end{lemma}

\proof
If $U \in \Inv{f}$ with $U \subset D$ then $\inv{U}{f} \subset D$, thus $\coinv{U}{f} \subset \co{D} = D$, which 
implies that $\coinv{U}{f} = D$, since $D$ is minimal. Hence $\coinv{U_1}{f} = \coinv{U_2}{f}$, and therefore by 
Lemmas \ref{lemma_5}~(7) and \ref{lemma_6}~(4) it follows that $U_1 \cap U_2 \ne \varnothing$.
\eop

We now start on the proof of Theorem~\ref{theorem_2}. Arbitrary intersections of $\bmod{f}$-invariant sets 
are $\bmod{f}$-invariant and $X$ is $\bmod{f}$-invariant and hence for each $A \subset X$ there is 
a least $\bmod{f}$-invariant set containing $A$, which will be denoted by $\pinv{A}{f}$. If $O$ is open then  by
Lemma~\ref{lemma_4} $\pinv{O}{f}$ is also open, since then $\interior{\pinv{O}{f}}$ is an 
$\bmod{f}$-invariant set containing $O$.

Suppose first $D \in \cofInv{f}$ is minimal. Since the topology on $X$ has a countable base then so does the 
subspace topology on $D$. There thus exists a sequence $\{O_m\}_{m\ge 1}$ of non-empty open subsets of  $D$ such 
that each non-empty open subset of $D$ contains some $O_m$. For each $m \ge 1$ put $V_m = \co{(\pinv{O_m}{f})}$; 
then $V_m \in \coInv{f}$ and $V \subset D$. Now for each $n \ge 1$ put $U_n = \bigcap_{m = 1}^n V_m$; then by 
Lemma~\ref{lemma_9} $U_n \ne \varnothing$ and so $U_n \in \coInv{f}$. Hence $\{U_n\}_{n \ge 1}$ is a 
decreasing sequence of sets from $\coInv{f}$ with $U_1 \subset D$. Consider $V \in \coInv{f}$ with 
$V \cap U_1 \ne \varnothing$; then $O_m \subset V \cap U_1$ for some $m \ge 1$ and so $V_m \subset V \cap U_1$, 
since $V \cap U_1 \in \coInv{f}$. In particular, $U_m \subset V$, since $U_m \subset V_m$, and this shows that 
$\gamma = \{U_n\}_{n \ge 1}$ is a transitive cascade. Moreover, $D \cap D_\gamma \supset U_1 \ne \varnothing$ and so 
$D \subset D_\gamma$. On the other hand, by Lemma~\ref{lemma_8} $D_\gamma = \coinv{U_1}{f}$, and 
$\coinv{U_1}{f} \subset D$, and hence $D_\gamma \subset D$. It follows that $D_\gamma = D$, which shows that each 
minimal element of $\cofInv{f}$ is the domain of a transitive cascade.

Conversely, let $\gamma = \{U_n\}_{n \ge 1}$ be a transitive cascade and by Lemma~\ref{lemma_7}
$D_\gamma \in \cofInv{f}$. Consider $E \in \cofInv{f}$ with $E \subset D_\gamma$. Then $U_m \subset E$ for some $m$, 
thus $\coinv{U_m}{f} \subset E$ and so by Lemma~\ref{lemma_8} $D_\gamma = \coinv{U_m}{f} \subset E$. This 
shows that $D_\gamma$ is minimal. 

This completes the proof of Theorem~\ref{theorem_2}.
\eop

If $U$ is a transitive element of $\coInv{f}$ then Theorem~\ref{theorem_2} implies in particular
that the domain $D_U$ of $U$ is a minimal element of $\cofInv{f}$. Moreover, by Lemma~\ref{lemma_8}
$D_U = \coinv{U}{f}$, and by Lemma~\ref{lemma_6}~(3) $\inv{U}{f} = \ninv{U}{f}$.

If $\gamma$ and $\gamma'$ are transitive cascades then by Theorem~\ref{theorem_2} the domains $D_\gamma$ 
and $D_{\gamma'}$, being minimal elements of $\cofInv{f}$, are either equal or disjoint. We next look at the 
relationship between $\gamma$ and $\gamma'$ when $D_\gamma = D_{\gamma'}$.

\begin{lemma}\label{lemma_10}
Let $\gamma = \{U_n\}_{n\ge 1}$ and $\gamma' = \{U'_n\}_{n\ge 1}$ be transitive cascades such that 
$D_\gamma = D_{\gamma'}$. Then the sequences $\gamma$ and $\gamma'$ are mutually cofinal, meaning that for each 
$p \ge 1$ there exist 
$q,\,r \ge 1$ with $U_q \subset U'_p$ and $U'_r \subset U_p$. In particular, we then have $C_\gamma = C_{\gamma'}$.
\end{lemma}

\proof 
Let $p \ge 1$. By Lemma~\ref{lemma_8} 
$\coinv{U_p}{f} \cap D_{\gamma'} = D_\gamma \cap D_{\gamma'} \ne \varnothing$ and hence by
Lemmas \ref{lemma_5}~(7) and \ref{lemma_6}~(4) $U_p \cap D_{\gamma'} \ne \varnothing$. Therefore 
$U'_r \subset U_p$ for some $r \ge 1$ and in the same way $U_q \subset U'_p$ for some $q \ge 1$, i.e., the 
sequences $\gamma$ and $\gamma'$ are mutually cofinal. 
\eop

Transitive cascades $\gamma$ and $\gamma'$ will be called \definition{equivalent} if the sequences $\gamma$ and 
$\gamma'$ are mutually cofinal. Lemma~\ref{lemma_10} implies that this is the case if and only if 
$D_\gamma = D_{\gamma'}$. By Theorem~\ref{theorem_2} each minimal element of $\cofInv{f}$ contains a unique 
equivalence class of transitive cascades. For each transitive cascade $\gamma = \{U_n\}_{n \ge 1}$ put
\[ 
\uninv{\gamma}{f} = \bigcap_{n\ge 1} \ninv{U_n}{f}\;.
\]
Therefore $x \in \uninv{\gamma}{f}$ if and only if for each $n \ge 1$ there exists an $m \ge 1$ such that
$f^m(x) \in U_n$. 
By Lemma~\ref{lemma_6}~(3) $\ninv{U_n}{f} = \inv{U_n}{f}$ for each $n$, and hence
$\uninv{\gamma}{f}$ is a fully $\bmod{f}$-invariant $G_\delta$-set with $\uninv{\gamma}{f} \subset D_\gamma$.
Moreover, by Lemma~\ref{lemma_8} $\coinv{U_n}{f} = D_\gamma$ for each $n \ge 1$, which implies that
\[ D_\gamma \setminus \uninv{\gamma}{f} = \bigcup_{n \ge 1} (\coinv{U_n}{f} \setminus \inv{U_n}{f})
\subset \bigcup_{n \ge 1} \partial \inv{U_n}{f}\;, \]
i.e., $D_\gamma \setminus \uninv{\gamma}{f}$ is contained in the meagre set $\bigcup_{n \ge 1} \partial \inv{U_n}{f}$.
If $U_m$ is transitive for some $m$ (and so $U_n = U_m$ for all $n \ge m$)
then $\uninv{\gamma}{f} = \ninv{U_m}{f}$, and in this case
$D_\gamma \setminus \uninv{\gamma}{f}$ is contained in the nowhere dense set $\partial \inv{U_m}{f}$.
If $\gamma$ and $\gamma'$ are equivalent then 
clearly $\uninv{\gamma}{f} = \uninv{\gamma'}{f}$; on the other hand, if $\gamma$ and $\gamma'$ are not equivalent 
then by Lemma~\ref{lemma_10} $\uninv{\gamma}{f} \cap \uninv{\gamma'}{f} = \varnothing$.

We now give a characterisation of the transitive elements in $\coInv{f}$ which corresponds to a standard result 
concerning the property of being topologically transitive (defined in terms of closed rather than open sets), and 
which can by found in Walters \cite{walters82}. For this we need an explicit expression for the least 
$\bmod{f}$-invariant set $\pinv{A}{f}$ containing $A$.

\begin{lemma}\label{lemma_11}
For each $A \subset X$
\[
\pinv{A}{f} = \{ x \in X : \mbox{$x = f^n(y)$ for some $y \in A$ and some $n \ge 0$} \}\;.
\]
\end{lemma}

\proof
Denote the set $\{ x \in X : \mbox{$x = f^n(y)$ for some $y \in A$ and some $n \ge 0$} \}$ by $A'$ for each 
$A \subset X$. We first show that $A'$ is $\bmod{f}$-invariant for each $A$. Let $x \in A' \setminus S$; then 
$f(x) \in X$ and $x \in f^n(A)$ for some $n \ge 0$, thus $f(x)\in X \cap f^{n+1}(A)$ and hence $f(x) \in A'$. 
Therefore $\bmod{f}(A') \subset A'$, i.e., $A'$ is $\bmod{f}$-invariant. Next note that if $B$ is 
$\bmod{f}$-invariant then by Lemma~\ref{lemma_2}~(2) $X \cap f^n(B) \subset B$ for all $n \ge 1$ and therefore
$B \subset B' = X \cap \bigcup_{n \ge 0} f^n(B) \subset B$; i.e., $B = B'$. Now let $A \subset X$ and $B$ be 
$\bmod{f}$-invariant with $A \subset B$; then $A' \subset B' = B$, and therefore $A' = \pinv{A}{f}$, since $A'$ is 
$\bmod{f}$-invariant and contains $A$.
\eop

In particular $\pinv{x}{f} = \{ y \in X : \mbox{$y = f^n(x)$ for some $n \ge 0$} \}$ for each $x \in X$. Thus if 
$x \in \Lambda_f(S)  = X \setminus \ninv{S}{f}$ then $\pinv{x}{f}$ is the set of points in the orbit of $x$ 
under $f$, and if $x \in \ninv{S}{f}$ then $\pinv{x}{f}$ is the finite set of points in the orbit of $x$ under 
$f$, but omitting the element $\bullet$.

\begin{proposition}\label{prop_3}
The following are equivalent for an element $U \in \coInv{f}$:

(1)\enskip  
$U$ is transitive.

(2)\enskip  $U \subset \co{(\ninv{O}{f})}$ for each non-empty open subset $O$ of $U$.

(3)\enskip  $U = \co{U}_*$, where $U_* = \{ x \in U : U = \co{(\pinv{x}{f})} \}$.

(4)\enskip  
$U = \co{(\pinv{x}{f})}$ for some $x \in U$.
\end{proposition}

\proof  
(1) $\Rightarrow$ (2):\enskip 
By Lemmas \ref{lemma_6}~(1) and \ref{lemma_4} the set
$X \setminus \ol{\ninv{O}{f}} = \interior{X \setminus \ninv{O}{f}}$ is $\bmod{f}$-invariant and so 
$V = U \setminus \ol{\ninv{O}{f}} = U \cap (X \setminus \ol{\ninv{O}{f}}\,)$ is either empty or an element of 
$\coInv{f}$. But $O \cap V = O \cap (U \setminus \ol{\ninv{O}{f}}) = O  \setminus \ol{\ninv{O}{f}} = \varnothing$
and it follows that $O \subset U \setminus \ol{V}$. In particular $U \ne V$, which implies $V = \varnothing$,
i.e., $U \subset \ol{\ninv{O}{f}}$. Thus $U \subset \co{(\ninv{O}{f})}$.

(2) $\Rightarrow$ (3):\enskip 
Since the topology on $X$ has a countable base so does the subspace topology on $U$. There thus exists a sequence 
$\{O_m\}_{m\ge 1}$ with each $O_m$ a non-empty open subset of  $U$ such that each non-empty open subset of $U$ 
contains some $O_m$. Let  $x \in U$; then  $U = \co{(\pinv{x}{f})}$ if and only if  
$\pinv{x}{f} \cap O_m \ne \varnothing$  for each  $m \ge 1$, and by Lemma~\ref{lemma_11} this holds if and 
only if $x \in \ninv{O_m}{f}$ for each  $m \ge 1$. This shows that $U_* = U \cap \bigcap_{m\ge 1}  \ninv{O_m}{f}$
(and in particular $U_*$ is a $G_\delta$-set). But for each $m \ge 1$ we have $U \subset \co{(\ninv{O_m}{f})}$ 
and hence
\[ 
U \cap \ninv{O_m}{f} = U \setminus (U \cap (\co{(\ninv{O_m}{f})} \setminus \ninv{O_m}{f}) 
\supset U \setminus \partial \ninv{O_m}{f}\;.
\]
It follows that 
$(X \setminus \ol{U}) \cup (U \cap \ninv{O_m}{f}) \supset X \setminus (\partial U \cup \partial \ninv{O_m}{f})$,
which implies that $(X \setminus \ol{U}) \cup (U \cap \ninv{O_m}{f})$ is a dense open subset of $X$. Therefore by 
the Baire category theorem the intersection
\[
(X \setminus \ol{U}) \cup U_* =   (X \setminus \ol{U}) \cup \Bigl(U \cap \bigcap_{m\ge 1}  \ninv{O_m}{f}\Bigr)
  = \bigcap_{m \ge 1} ((X \setminus \ol{U}) \cup (U \cap \ninv{O_m}{f})) 
\]
is also dense and so
\[ 
X = \ol{(X \setminus \ol{U}) \cup U_*} = \ol{(X \setminus \ol{U})} \cup \ol{U_*} 
= (X \setminus \co{U}) \cup \ol{U_*} = (X \setminus U) \cup \ol{U_*}\;.
\]
This shows that $U \subset \ol{U_*}$, i.e., $U \subset \co{U}_*$.

(3) $\Rightarrow$ (4):\enskip This is clear, since $U \cap U_* \ne \varnothing$.

(4) $\Rightarrow$ (1):\enskip 
Let $V \in \coInv{f}$ with $V \subset U$ and let $x \in U$  be such that $U = \co{(\pinv{x}{f})}$. There thus 
exists  $m \ge 0$  with  $f^m(x) \in V$ and then $\pinv{f^m(x)}{f} \subset V$, since  $V$ is $\bmod{f}$-invariant. 
Put $B = \{x,f(x),\ldots,f^{m-1}(x)\}$; then
\[ 
U = \co{(\pinv{x}{f})} \subset \ol{B \cup \pinv{f^m(x)}{f}}  = B \cup \ol{\pinv{f^m(x)}{f}} \subset B \cup  \ol{V}\;,
\]
and so  $U \setminus B' \subset \ol{V}$ with $B' = B \setminus \ol{V}$. It follows that $U \subset V$, since $B'$ 
is finite and $X$ is perfect. (If $W$ is open with $W \cap U \ne \varnothing$ then 
$W \cap (U \setminus B') \ne \varnothing$ and so $W \cap V \ne \varnothing$, since $U \setminus B' \subset \ol{V}$. 
Hence $U \subset \interior{\ol{V}} = \co{V} = V$.)
\eop

We now turn to another decomposition of $X$ into disjoint elements of $\fInv{f}$ which more directly involves the 
set $S$. Each subset $C$ of $S$ is $\bmod{f}$-invariant, since $\bmod{f}(C) = f(C \setminus S) = \varnothing$, 
and therefore by Lemma~\ref{lemma_6}~(3) $\inv{C}{f} = \ninv{C}{f}$, i.e.,
\[  \inv{C}{f} = \{ x \in X : \mbox{$f^n(x) \in C$ for some $n \ge 0$} \}\;. \]
In particular, $X$ is the disjoint union of the sets $\inv{S}{f}$ and $\Lambda_f(S)$, where as before
$\Lambda_f(S) = \{ x \in X : \mbox{$f^n(x) \notin S$ for all $n \ge 0$} \}$.

Put $\Zed{f} = \interior{\Lambda_f(S)}$ and $\Sigma_f = \coinv{S}{f}$. Then we have $\Zed{f} \in \fInv{f}$,
$\Sigma_f$ is either empty or an element of $\cofInv{f}$, the sets $\Zed{f}$ and $\Sigma_f$ are disjoint
and their union $\Zed{f} \cup \Sigma_f$ is a dense open subset of $X$, since
$X \setminus (\Zed{f} \cup \Sigma_f) = \partial\, \ol{\inv{S}{f}}$ and the boundary of a closed set is nowhere
dense. If $D$ is a minimal element of $\cofInv{f}$ then either $D \subset \Sigma_f$ or 
$D \cap \Sigma_f = \varnothing$ (in which case $D \subset \co{Z}_f$).

\begin{theorem}\label{theorem_3}
Suppose $S$ is finite and $\Sigma_f \ne \varnothing$ (and so $S$ is non-empty). Then $\Sigma_f$ contains at least 
one and at most $\#(S)$ minimal elements of $\cofInv{f}$ (with $\#(S)$ the cardinality of $S$).
\end{theorem}

\proof
This requires some preparation.

\begin{lemma}\label{lemma_12}
If $C \subset S$ and $U \in \Inv{f}$ with $U \cap \coinv{C}{f} \ne \varnothing$ then $U \cap C \ne \varnothing$. 
In particular, if $\coinv{C}{f} \ne \varnothing$ then $\coinv{C}{f}$ contains an element of $C$, and each 
$U \in \coInv{f}$ with $U \subset \Sigma_f$ contains an element of $S$.
\end{lemma}

\proof 
Put $V = U \cap \coinv{C}{f}$. Then by Lemma~\ref{lemma_5}~(7) $V \cap \inv{C}{f} \ne \varnothing$. Thus, 
since $\inv{C}{f} = \ninv{C}{f}$, there exists $c \in C$, $x \in V$ and $n \ge 0$ with $f^n(x) = c$ and it 
follows that $c = f^n(x) \in X \cap f^n(V)$. But by Lemma~\ref{lemma_3} $X \cap f^n(V)$ is open and 
$X \cap f^n(V) \subset X \cap f^n(U) \subset U$. Hence $c \in U$, which shows $U \cap C \ne \varnothing$.
\eop

\begin{proposition}\label{prop_4}
If $c \in S$ with $\coinv{c}{f} \ne \varnothing$ then $\coinv{c}{f}$ is a minimal element of $\cofInv{f}$ 
containing $c$.
\end{proposition}

\proof 
By Lemma~\ref{lemma_12} $c \in \coinv{c}{f}$.
Let $D \in \cofInv{f}$ with $D \subset \coinv{c}{f}$. Then, again by Lemma~\ref{lemma_12}
$c \in D$, hence $\inv{c}{f} \subset D$, since $\inv{c}{f}$ is the least fully $\bmod{f}$-invariant
set containing $c$, and then $\coinv{c}{f} \subset \co{D} = D$. Thus $D = \coinv{c}{f}$, which shows that
$\coinv{c}{f}$ is minimal. 
\eop

\begin{lemma}\label{lemma_13}
Let $C_1,\,\ldots,\,C_n$ be subsets of $S$ and put $C = \bigcup_{k=1}^n C_k$. Then
\[ \bigcup_{k = 1}^n \coinv{C_k}{f} \subset \coinv{C}{f} \subset
     \bigcup_{k = 1}^n \coinv{C_k}{f} \cup N \]
with $N$ the nowhere dense set $\,\bigcup_{k=1}^n \partial\, \ol{\inv{C_k}{f}}$.
\end{lemma}

\proof
Since $C_k \subset C$ it immediately follows that $\coinv{C_k}{f} \subset \coinv{C}{f}$ and hence that
$\bigcup_{k = 1}^n \coinv{C_k}{f} \subset \coinv{C}{f}$. On the other hand
$\inv{C}{f} = \bigcup_{k=1}^n \inv{C_k}{f}$ and thus
\begin{eqnarray*}
\coinv{C}{f} \subset \ol{\inv{C}{f}} &=& \bigcup_{k=1}^n \ol{\inv{C_k}{f}} 
= \bigcup_{k=1}^n \bigl(\coinv{C_k}{f} \cup \bigl(\,\ol{\inv{C_k}{f}} \setminus \coinv{C_k}{f}\bigr)\bigr)\\ 
&=& \bigcup_{k=1}^n \bigl(\coinv{C_k}{f} \cup \partial\,\ol{\inv{C_k}{f}}\bigr)  
= \bigcup_{k=1}^n \coinv{C_k}{f} \cup N\;.\ \eop
\end{eqnarray*}

We now come to the proof of Theorem~\ref{theorem_3}. First,
since $S = \bigcup_{c \in S}$ it follows from Lemma~\ref{lemma_13} that
$\bigcup_{c \in S} \coinv{s}{f} \subset \Sigma_f \subset \bigcup_{c \in S} \coinv{c}{f} \cup N$, where $N$ is nowhere
dense. Thus, since $\Sigma_f$ is non-empty and open, there exists a $c \in S$ such that 
$\coinv{c}{f} \ne \varnothing$, and then by Proposition~\ref{prop_4} $\coinv{c}{f}$ is a minimal element 
of $\cofInv{f}$, which shows that $\Sigma_f$ contains at least one minimal element of $\cofInv{f}$. Now let $D$ be 
any minimal element of $\cofInv{f}$ contained in $\Sigma_f$. Then $D \cap \coinv{c}{f} \ne \varnothing$ for some 
$c \in S$ (since $N$ is nowhere dense) and then $D = \coinv{c}{f}$, since by Proposition~\ref{prop_4} 
$\coinv{c}{f}$ is also a minimal element of $\cofInv{f}$. Therefore $\Sigma_f$ contains at most $\#(S)$ minimal 
elements of $\cofInv{f}$. This completes the proof of Theorem~\ref{theorem_3}.
\eop

Theorems \ref{theorem_2} and \ref{theorem_3}
together show that if $S$ is finite and $\Sigma_f \ne \varnothing$ then $\Sigma_f$ contains at least 
one and at most $\#(S)$ equivalence classes of minimal cascades.

\begin{proposition}\label{prop_5}
If $D \in \cofInv{f}$ with $D \subset \Sigma_f$ then $D = \coinv{S \cap D}{f}$ for each $D \in \cofInv{f}$.
\end{proposition}

\proof 
First consider an element $U \in \coInv{f}$ with $U \subset \Sigma_f$. Then, since $S \cap U \subset U$, it 
follows that $\inv{S \cap U}{f} \subset \inv{U}{f}$ and therefore that $\coinv{S \cap U}{f} \subset \coinv{U}{f}$. 
Now if the $\bmod{f}$-invariant open set $V = U \setminus \ol{\inv{S \cap U}{f}}$ were non-empty then by 
Lemma~\ref{lemma_12} it would contain an element of $S \cap U$. Since $S \cap U \subset \ol{\inv{S \cap U}{f}}$ 
this is not possible, and so $V = \varnothing$, i.e., $U \subset \ol{\inv{S \cap U}{f}}$. Hence 
$\inv{U}{f} \subset \ol{\inv{S \cap U}{f}}$, since $\ol{\inv{S \cap U}{f}}$ is fully $\bmod{f}$-invariant, and 
thus also $\inv{U}{f} \subset \coinv{S \cap U}{f}$. This shows $\coinv{U}{f} = \coinv{S \cap U}{f}$. In particular,
if $D \in \cofInv{f}$ with $D \subset \Sigma_f$ then $D = \coinv{S \cap D}{f}$, since here $D = \coinv{D}{f}$.
\eop

By Proposition~\ref{prop_5} the mapping $D \mapsto S \cap D$ from the set 
$\{ D \in \cofInv{f} : D \subset \Sigma_f \}$ to the set of subsets of $S$ is injective.

\bigskip
\bigskip

{\sc Fakult\"at f\"ur Mathematik, Universit\"at Bielefeld}\\
{\sc Postfach 100131, 33501 Bielefeld, Germany}\\
\textit{E-mail address:} \texttt{preston@math.uni-bielefeld.de}\\

\end{document}